\documentclass[letterpaper, 10 pt, conference]{ieeeconf}
\usepackage{cite}
\usepackage{amsmath,amssymb,amsfonts}
\usepackage{graphicx}
\usepackage{hyperref}
\usepackage{textcomp}
\usepackage{cleveref}
\usepackage{mathrsfs}
\usepackage{bm}
\usepackage[normalem]{ulem}
\usepackage{subfiles}
\usepackage{bbm}
\usepackage{comment}
\let\labelindent\relax
\usepackage[shortlabels]{enumitem}
\usepackage{xcolor}
\usepackage{hyperref}
\usepackage[normalem]{ulem}
\usepackage{tikz}
\usetikzlibrary{patterns}
\PassOptionsToPackage{hyphens}{url}
\usepackage{url}
\usepackage{caption}

\usepackage{algorithm}
\usepackage{algpseudocode}

\definecolor{ao}{rgb}{0.0, 0.5, 0.0}
\usepackage{todonotes}

\newtheorem{assumption}{Assumption}

\usepackage{times}

\DeclareMathOperator*{\argmin}{arg\,min}
\newcommand{\norm}[1]{\left\lVert#1\right\rVert}
\renewcommand{\vec}[1]{{\boldsymbol{#1}}}

\newcommand{\R}{\mathbb{R}}
\newcommand{\E}{\mathbb{E}}
\newcommand{\N}{\mathbb{N}}

\newcommand{\cF}{{\mathcal F}}

\newcommand{\cL}{{\mathcal L}}

\newcommand{\cP}{{\mathcal P}}

\newcommand{\cS}{{\mathcal S}}

\newcommand{\cX}{{\mathcal X}}

\newcommand{\Nmax}{\Delta}

\newcommand{\paran}[1]{\left(#1\right)}
\newcommand{\inp}[2]{\langle #1, #2 \rangle}
\newcommand{\red}[1]{\textcolor{red}{#1}}
\newcommand{\fopt}{f^{\star}}
\newcommand{\powercoeff}[1]{\sigma_{#1}}

\newcommand{\snr}[1]{\rho_{#1}}
\renewcommand{\red}[1]{\textcolor{black}{#1}}

\definecolor{OliveGreen}{rgb}{0,0.5,0}
\newcommand{\av}[1]{{\color{black} #1}}

\newcommand{\del}[2]{\frac{\partial #1}{\partial #2}}
\newcommand{\Pavg}{\bar\sigma_{p}^2}
\newcommand{\Pavgcoeff}{\bar\sigma_{p}}

\newcommand{\Cmin}{\powercoeff{\text{lb}}}

\def \enc {\texttt{enc}}
\def \dec {\texttt{dec}}
\newcommand{\contr}{\gamma_{\mu\eta}}
\newcommand{\tswitch} {t_{\text{switch}}}
\newtheorem{theorem}{Theorem}

\newtheorem{corollary}{Corollary}

\newtheorem{definition}{Definition}

\newtheorem{lemma}{Lemma}

\newtheorem{problem}[theorem]{Problem}
\newtheorem{proposition}[theorem]{Proposition}
\newtheorem{remark}{Remark}

\def\BibTeX{{\rm B\kern-.05em{\sc i\kern-.025em b}\kern-.08em
    T\kern-.1667em\lower.7ex\hbox{E}\kern-.125emX}}
\title{Power-Constrained Policy Gradient Methods for LQR}
\author{Ashwin Verma, Aritra Mitra, Lintao Ye, and Vijay Gupta
\thanks{A. Verma and V. Gupta are with the Purdue University (e-mail: \{verma240, gupta869\}@purdue.edu).  A. Mitra is 
with North Carolina State University (e-mail: amitra2@ncsu.edu). L. Ye is with Huazhong University of Science and Technology (e-mail: yelintao93@hust.edu.cn)}}
\IEEEoverridecommandlockouts
\begin{document}
\maketitle

\begin{abstract}
Consider a discrete-time Linear Quadratic Regulator (LQR) problem solved using policy gradient descent when the system matrices are unknown. The gradient is transmitted across a noisy channel over a finite time horizon using analog communication by a transmitter with an average power constraint. This is a simple setup at the intersection of reinforcement learning and networked control systems. We first consider a communication-constrained optimization framework, where gradient descent is applied to optimize a non-convex function under noisy gradient transmission. We provide an optimal power allocation algorithm that minimizes an upper bound on the expected optimality error at the final iteration and show that adaptive power allocation can lead to better convergence rate as compared to standard gradient descent with uniform power distribution. We then apply our results to the LQR setting.
\end{abstract}

\section{Introduction}

There is a recent surge of interest in model-free approaches to the Linear Quadratic Regulator (LQR) problem. Among such methods, policy gradient (PG) algorithms, in particular, have gained significant popularity due to their simplicity and practical applicability. When applied to the classical LQR problem \cite{anderson}, prior work \cite{fazel} has shown that despite the inherent non-convexity of the optimization landscape, PG algorithms with noise-free gradient estimates can guarantee convergence to the globally optimal policy. Here, we are interested in the utility of such algorithms to the setup shown in networked control systems. Specifically, we would like to characterize the robustness of these algorithms to communication-induced distortions, which arise when gradient or policy updates are transmitted over realistic communication channels. For gradients transmitted over noise-free but quantized channels, \cite{pmlr-v242-mitra24a} establishes a somewhat surprising result that when the bit rate exceeds a certain threshold, there exist algorithms that ensure exponentially fast convergence to the optimal policy, with {\em no degradation} in the convergence rate compared to the unquantized setting. In this work, we investigate a similar problem for noisy analog channels with average power constraints. 

We consider a setup in which the policy gradients computed by a worker agent are transmitted to the decision maker (or a server) over a noisy channel, who then updates the policy. Our goal is to design $(i)$ a power allocation scheme at the worker and $(ii)$ a policy update rule at the decision-maker, so that the resulting policy gradient algorithm continues to guarantee convergence to the neighborhood of optimal solution. Specifically, we are interested in the application of the gradient descent method to a learning problem for LQR under communication constraints on policy gradient updates. 

The problem is closely related to a communication-constrained optimization problem in which gradient descent is applied to optimize a non-convex function under a similar noisy gradient transmission. We also begin by analyzing such a setup. Gradient-based methods are widely used in optimization and control due to their simplicity, computational and memory efficiency, and robustness  \cite{gasnikov2017universal,belkin2021fit,nesterov2015universal,karimi2016linear}. In keeping with our final goal of the LQR problem, we focus on minimizing non-convex functions that satisfy the well-known Polyak–{\L}ojasiewicz (PL) condition \cite{polyak1963gradient} and have Lipschitz continuous gradients. In particular, for the LQR problem, a well known challenge is that the objective function defined over the policy space, has gradients that satisfy the PL condition and smoothness properties only \emph{locally}. 

We note two related lines of work here. The first direction is in works such as~\cite{konevcny2016federated} that study gradient-based optimization in which a central parameter server executing the gradient iteration has access only to noisy gradient estimates that have been transmitted over a communication channel by an oracle or worker agent. The second direction studies such methods when communication involves over-the-air transmission subject to certain power constraints\cite{jha2022fundamental,zhang2020gradient, liang2024communication}. However, most of these existing works consider per-iteration power constraints, limiting the transmission power at each gradient update. In contrast, we investigate a scenario in which the communication is subject to an average power constraint over the entire optimization process, \av{similar to the setting in \cite{sun2021dynamic}.} 

This formulation allows for power accumulation, enabling more refined gradient transmissions at critical iterations. We note that for functions satisfying {\em global} properties,~\cite{sun2021dynamic} considered average power constraint for the problem of federated learning with over-the-air communication but the focus of that work was the design of a dynamic device scheduling algorithm. We also note that unlike these works that consider a standard optimization setup where the primary concern is the final iterate value, our primary motivation is the LQR problem which introduces additional challenges. Besides the fact that in LQR, PL, and Lipschitz properties are satisfied only locally, it is also crucial to control how the updates evolve throughout the process. This necessitates a more careful algorithm to ensure stable and efficient learning under power constraints. 

\textbf{Outline and Contributions:} 
In Section~\ref{sec:prob_form}, we formally state the \av{LQR problem} considered and describe the policy gradient algorithm. In Section~\ref{sec:opt_prob}, we provide the problem formulation in an equivalent optimization landscape.
 \av{In optimization landscape we consider the problem of power allocation for Power-Allocated Gradient Descent (\texttt{PAGD}) of functions satisfying regularity conditions globally and locally.} 

In Sections~\ref{subsec:opt_pow_all} and ~\ref{subsec:opt_pow_all2}we state two main results, Theorems~\ref{thm:PL} and \ref{thm:local_conv}, regarding optimal power allocation. 
Theorem~\ref{thm:PL} provides the convergence result for functions that satisfy the PL and Lipschitz properties globally. 
Theorem~\ref{thm:local_conv} provides the analysis for functions that satisfy the desired properties only within a local neighborhood. To ensure that the updates remain in this region, we introduce an additional constraint on the power allocation of the form $\powercoeff{t}\geq \powercoeff{lb}$ and provide a lower bound for $\powercoeff{lb}$. 
\av{Finally, we state the allocation scheme and convergence results for the LQR problem in Theorem~\ref{thm:lqr_thm}.}

\textbf{Notation:} For any positive integer $n$, let $[n] = \{1,2,\dots,n\}$ and $[n]_0 = [n]\cup \{0\}$. Denote the set of all positive definite $n\times n$ matrices by $\mathbb{S}^{n}_{++}$. For vector $x \in \R^n$, denote its Euclidean norm by $\|x\|$. For a matrix $A$, we use the same notation $\|A\|$ to denote its Frobenius norm; the distinction will be clear from context. We denote the inner product between vectors $x,y \in R^n$ by $\inp{x}{y}$. The notation $\{z_t\}_{t\ge0}$ defines a sequence $z_{t}$ over the times $t=0,1,\cdots.$

\section{Problem Formulation:}\label{sec:prob_form}

Consider a remote sensing agent that transmits analog data to a decision-maker across a communication channel that adds noise to any transmitted signal. The agent is allocated a limited average power budget and needs to allocate the power to the signal sent at each transmission. 

{\textbf{The LQR problem:}} Specifically, consider the linear time-invariant (LTI) 
\[x_{t+1} = Ax_t + Bu_t + w_t,\qquad t\geq 0,\]
where $A\in \R^{n\times n}$ and $B\in \R^{n\times m}$ are system matrices, $x_t$ and $u_t$ are the state and control input vectors at time $t$, and $\{w_t\}_{t\ge0}$ is a random process with independent and identically distributed random variables
. The pair $(A,B)$ is assumed to be controllable. 
Without loss of generality, we assume that the initial state is $x_0 = 0$. 
The LQR problem aims to design  $\{u_t\}$ that minimizes the average cost function,  
$$
\lim_{T\to \infty} \frac{1}{T}\E[\sum_{t=0}^{T-1}x_t^T Q x_t + u_t^T Ru_t],
$$
where $Q \in \mathbb{S}^n_{++}$ and $R \in \mathbb{S}_{++}^m$ are cost matrices and expectation is taken with respect to the disturbance process $\{w_t\}_{t\ge 0}$. It is well-known~\cite{bertsekas2015dynamic} that the optimal control inputs are given by the static state-feedback policy $u_t = Kx_t$ for a stabilizing controller $K\in \R^{n\times m}$ as given by 
\begin{align}\label{eqn:LQR obj J(K)}
\!K^{\star} \!\!=\!\! \argmin_{K} J(K) = \argmin_{K} \text{trace}((Q+K^{\top}RK)\Sigma_K),    
\end{align}

where $\Sigma_K\in\mathbb{S}_{++}^n$ is the solution to the Riccati equation: $$\Sigma_K = \Sigma_w + (A+BK)^{\top} \Sigma_K (A+BK).$$

\textbf{Policy gradient to solve LQR:} When the system matrices $A$ and $B$ are unknown, the optimal $K^{\star}$ can be obtained through a policy gradient method, 
$K_{t+1} = K_t - \eta \nabla J(K_t)$, initialized with an arbitrary stabilizing matrix $K_0$ and a well-chosen step-size $\eta > 0$ \cite{fazel,malik2020derivative}. With known matrices $A, B, Q$, and $R$, the exact gradient $\nabla J(K_t)$ can be computed and the algorithm converges exponentially fast to the optimal policy $K^{\star}$. 
When the matrices $A$ and $ B $ are unknown accurate estimates of $J(K)$ and $\nabla J(K)$ can be computed through the system trajectories obtained by applying the control policy $u_t = K x_t$. 

\textbf{Communication Constraint:}
At each time $t\in[0,T]$, the agent determines the gradient ${\nabla J(K_t)}$ of the function at the current value of the state variable and transmits this gradient to the decision-maker across a noisy communication channel. 
We denote the signal transmitted by the agent as $\enc(\nabla J(K_t))$ to reflect the fact that it is an encoding of the gradient. We assume that the agent utilizes analog modulation to communicate gradient information. Further, in anticipation of the fact that we will impose a power constraint at the transmitter, we note that without loss of generality and to save power, the transmitted signal can be first normalized by $G$. Finally, we assume that the transmitter must satisfy an average power constraint. Thus, if the power allocated at time step $t$ is denoted by $\sigma_{t}^{2},$ then the encoded transmission is given by
$$\enc(\nabla J(K_t)) = \powercoeff{t}\frac{\nabla J(K_t)}{G},$$
and the average transmission power must satisfy
    \begin{align}\label{eqn:avg_power_constraint}
        \frac{\sum_{t=0}^{T-1}\powercoeff{t}^2}{T} \leq \Pavg.
    \end{align}
The communication channel adds a noise term to the transmitted signal, so that the received signal at the output of the channel is given by 
 $$g_t = \enc(\nabla(J(K_t)))+n_t.$$   
 We assume the following. 
\begin{assumption}\label{asmp:noise}
The stochastic sequence $\{n_t\}$ consists of independent random variables that have zero mean, ${\E[n_t] = \vec{0}}$, and bounded second moment for the norm, $\E(\norm{n_t}^2) = \sigma_N^2$, any $t\in \N_0$.

\end{assumption}
Furthermore, to ensure that the updates remain within the set of stabilizing matrices, we impose the following assumption on the noise.
\begin{assumption}\label{asmp:noise_bounded}
For all $t \in [T-1]_0$, the random variables $n_t$ satisfies almost-sure boundedness, ${\Pr(\norm{n_t} \leq \Nmax) = 1}$.
\end{assumption}
The decision maker receives the signal $g_t$ and decodes it to get an estimate of the gradient $\dec(g_t)$ as $$\dec(g_t) = \frac{G}{\powercoeff{t}}g_t.$$
Using this value, it performs a gradient-descent step with step-size $\eta$ to update the state variable as 
$$K_{t+1} = K_t - \eta \, \dec(g_t).$$ 
The decision-maker then sends the updated variable, $K_{t+1}$, to the agent over a noiseless channel. The noiseless assumption for the transmission by the decision-maker is justified by the fact that it is a more resourceful agent with sufficient power. 
\av{We refer to the gradient descent algorithm as Power-Allocated Gradient Descent (\texttt{PAGD}). The \texttt{PAGD} is characterized by the power allocation scheme $\{\powercoeff{t}\}_{t=0}^{T-1}. $}

\textbf{Problem considered:} 
\av{\begin{problem}[$\cP_0$]
For the LQR problem~\eqref{eqn:LQR obj J(K)} given the total power budget $T\Pavg$, determine a power allocation scheme $\{\powercoeff{t}\}_{t=0}^{T-1}$ for \texttt{PAGD} to minimize $\E[J(K_T) - J^{\star}]$. 
\end{problem}

Since directly minimizing $\E[J(K_T)-J^{\star}]$ depends on the full knowledge of the problem parameters, we instead focus on minimizing a tractable upper bound on the error. This bound depends on power allocated $T\Pavg$ and certain properties, such as, upper bound on the maximum singular value of matrices $A, B$, of the problem parameters and avoids needing the entire matrix.}

\section{Optimization Problem}\label{sec:opt_prob}

\textbf{Optimization setup:} As a step towards solving the problem $\mathcal{P}_{0}$, we first pose and solve a noisy power-constrained optimization problem. In this problem our objective is to minimize a function ${f:\R^d \to \R}$ using a gradient algorithm implemented from $t=0,\cdots,T-1$ for a given horizon $T$. 
\av{Note that even though $K$ and $\nabla J(K)$ are matrices the analysis of the optimization holds true since we can vectorize the matrices and apply the descent algorithm.
For the optimization problem, we begin by assuming that the function $f$ satisfies the following properties over its entire domain. }
 \begin{assumption}[Smoothness]
 \label{assumption_smooth}
  The function $f$ is continuously differentiable. Further, it is $L$-smooth for a constant $L>0$, so that the gradient map $\nabla f : \R^d \to \R^d$ is $L$-Lipschitz continuous, i.e., ${\|\nabla f(x) - \nabla f(y)\| \leq L \|x - y\|}, \quad \forall x, y \in \R^d.$
\end{assumption}
\begin{assumption}[Polyak–{\L}ojasiewicz (PL) Condition]
 \label{assumption_PL}
The function $f$ satisfies the $\mu$-PL condition over the domain $\cX$ for some constant $\mu > 0,$  so that 
    \begin{align}
        f(x) - f^\star \leq \frac {1}{2\mu} \| \nabla f(x) \|^2, \quad \forall x \in \cX,
    \end{align}
    where $f^\star = f(x^{\star})$ is the function value at the global solution $x^{\star}$. 
\end{assumption}
\begin{assumption}[Bounded Gradients]\label{asmp:bounded_gradient}
    $f$ has uniformly bounded gradients over the domain $\cX$, i.e., $\|\nabla f(x)\| \leq G $.
\end{assumption} 
\av{We follow similar definitions for the encoder, decoder, and the noise sequence as set for problem $\cP_0$. 
Specifically, at each time $t\in[0,T]$, an agent determines the gradient ${\nabla f(x_t)}$ of the function at the current value of the state variable and transmits this gradient to a decision-maker across a noisy communication channel as $\enc(\nabla f(x_t))$.

We assume that the transmitter must satisfy an average power constraint. Thus, if the power allocated at time step $t$ is denoted by $\sigma_{t}^{2},$ then the encoded transmission is given by
$\enc(\nabla f(x_t)) = \powercoeff{t}\frac{\nabla f(x_t)}{G},$
and the average transmission power must satisfy eq.~\eqref{eqn:avg_power_constraint}.

As in $\cP_0$, the communication channel adds a noise term to the transmitted signal, so that the received signal at the output of the channel is given by 
$g_t = \enc(\nabla(f(x_t)))+n_t.$
The decision maker decodes $g_t$ to get an estimate of the gradient $\dec(g_t)$ as $\dec(g_t) = \frac{G}{\powercoeff{t}}g_t.$
Using this value, it performs a gradient-descent step with step-size $\eta$ to update the state variable as 
$x_{t+1} = x_t - \eta \dec(g_t).$
The decision-maker then sends the updated variable, $x_{t+1}$, to the agent over a noiseless channel. The noise sequence $\{n_t\}_{t\geq0}$ follows assumption~\ref{asmp:noise}. We refer to the descent algorithm also as \texttt{PAGD} algorithm.}

\av{
\textbf{Proposed Allocation: } In optimizing the power allocation for a bound on the last-iterate error, we derive a power allocation structure which we refer to as \emph{Constant-then-Geometric} (CtG) power allocation, as described in Allocation Scheme~\ref{alg:poweralloc}. This approach involves assigning a constant, baseline power level during the initial phase of the algorithm, $\forall t \in [\tswitch-1]_0$, followed by a geometrically increasing power allocation in the latter stages.

The CtG allocation balances robustness and precision by initially using constant power to ensure stability and prevent divergence when iterates are far from optimal. This phase conserves energy and maintains a minimum SNR, which is especially useful when only local convergence guarantees are available. 

As the iterates approach the optimum, precision becomes critical. Power is then increased geometrically, effectively improving gradient accuracy without reducing the step-size. This mirrors the benefits of step-size decay via enhanced communication quality.
When the function satisfies global smoothness and PL conditions, the optimal allocation reduces to a fully geometric scheme, corresponding to CtG allocation with $\tswitch = 0$. In contrast, for locally constrained functions, the optimal switch time from constant to geometric allocation is determined based on problem parameters $\mu$, $L$, the total power budget $T\Pavg$, and algorithmic hyperparameters $\eta$ and $\powercoeff{lb}$.
CtG allocation thus offers a principled, non-adaptive allocation scheme that minimizes a standard upper bound on expected error and performs effectively under both global and local constraints.
}

\subsection{Global Constraints}\label{subsec:opt_pow_all}

\av{\begin{problem}[$\cP_1$]
Consider a function $f:\R^d\to \R$ satisfying assumptions~\ref{assumption_smooth}, \ref {assumption_PL}, and \ref{asmp:bounded_gradient} over $\R^d$.
Given the total power budget $T\Pavg$ determine a power allocation scheme $\{\powercoeff{t}\}_{t=0}^{T-1}$ for \texttt{PAGD} to minimize $\E[f(x_T) - \fopt]$.
\end{problem}}
\av{Similar to the argument for $\cP_0$, we focus on minimizing a tractable upper bound on the suboptimality error to determine a power allocation scheme}
The problem $\mathcal{P}_{1}$ is to design the power  $\sigma^{2}_{t}$ allocated at each time $t$ in a way that satisfies the constraint~(\ref{eqn:avg_power_constraint}) and the expected last-iterate error $\E[f(x_T) - \fopt]$ is minimized. 
\av{
\begin{remark}
    Although we assume gradient boundedness over the entire domain $\mathbb{R}^d$, as done in \cite{koloskova2019decentralized}, our analysis only requires this condition to hold along the optimization trajectory, that is, $\|\nabla f(x_t)\| \leq G$ for all $t\in [T-1]_0$.
\end{remark}
}

To gain some intuition into the problem, we consider the case when no optimization of the allocated power is done. Recall that the gradient descent iterates are performed as
\begin{align}\label{eq:updateRule}
    {x_{t+1} =  x_t - \eta \paran{ {\nabla f(x_t)} + \frac{G}{\powercoeff{t}}n_t }.}
\end{align}
In the absence of any constraints on the average power that can be allocated, we approach the classical (noiseless) gradient. The simplest approach for the allocation of power is to utilize constant power for each transmission, i.e., $\powercoeff{t} = \Pavgcoeff$ for all $t\in \N_0$. In this case, we can express the gradient  descent iterates as $$x_{t+1} = x_t - \eta(\nabla f(x_t) + n_t'),$$ where $\{n_t'\}_{t\geq 0}$ is a sequence of random variables that satisfies Assumption~\ref{asmp:noise} with variance $\sigma_{N'}^2 = \frac{G^2\sigma_N^2}{\Pavg}$. \av{Following the results for SGD such as \cite[Theorem 4.6]{gower2021sgd} we get the following bound on the expected error for the last-iterate of the algorithm. }
\begin{proposition}{\av{(Following \cite[Theorem 4.6]{gower2021sgd}})}\label{thm:constant_power} 
    Consider Problem $\mathcal{P}_{1}$ specified above. Let ${\eta \in \left(0,\frac{1}{L}\right) }$ and the power allocation $\powercoeff{t} = \Pavgcoeff$ for all $t\in [T-1]_0$. 
    The last-iterate of the gradient descent satisfies 
    \begin{align*}
        \E[f(x_T) - \fopt] \leq \paran{1-\mu\eta}^T (f(x_0)-\fopt) + \frac{LG^2\eta}{\mu }\frac{\sigma_N^2}{\Pavg}.
    \end{align*}
\end{proposition}
Note that there are two components to the upper bound on the expected error -- the exponentially decaying error of the initial estimate error and the constant error due to the presence of noise. Since we utilize a constant step-size, the gradient descent algorithm will achieve only limited accuracy, leading the function value at the iterates to a neighborhood of the optimal point. 

\begin{algorithm}[t!]
\caption{Power-Allocated Gradient Descent (\texttt{PAGD})}
\label{alg:PAGD}
\begin{algorithmic}[1]
\State \textbf{Initialization:} $x_0 = 0$.
\State $\{\powercoeff{t}\}_{t=0}^{T-1} = \texttt{CtG}(T, \Pavg, \mu, \eta, \powercoeff{lb})$
\State \textbf{for} \hspace{0.25mm} {$t \in [T-1]_0$} \hspace{0.25mm} \textbf{do}
\State \quad \textbf{At Worker:}
 \State \quad Receive iterate $x_t$ and gradient $g_{t-1}$ from server. 
 \State \quad Compute $\nabla f(x_t)$ and transmit $\powercoeff{t} \frac{\nabla f(x_t)}{G}$.
 \State \quad \textbf{At Decision-Maker/Server:}
 \State \quad Receive $g_t = \powercoeff{t} \frac{\nabla f(x_t)}{G} + n_t$.
 \State \quad Update the model as: 
 \begin{align}
     x_{t+1} = x_t - \eta \, \dec({g_t}).
 \end{align}
\State \textbf{end for}
\State \Return $x_T$
\end{algorithmic}
\end{algorithm}

\av{\begin{algorithm}[t!]
\captionsetup{name=Allocation}
\caption{\texttt{CtG}$(T, \Pavg, \mu, \eta, \powercoeff{lb})$}
\begin{algorithmic}[1]
\State \textbf{Known Parameters:} $\mu, \Pavg, T$
\State \textbf{Hyperparameters:} $\eta, \powercoeff{lb}$
\State $\contr = \sqrt{1- \mu \eta}$, 
\State \textbf{if} $\Pavg < \powercoeff{lb}^2$
\State \quad \Return Error: Insufficient budget for power allocation
\State ${\tswitch = \min\{t \in [T-1]_0 : \contr^{T-1-t}\geq \frac{1 - \contr^{T-t}}{1-\contr} \frac{\powercoeff{lb}^2}{T\Pavg - t \powercoeff{lb}^2}\}}$
\State \textbf{for} $t \in [T-1]_0$ \textbf{do}
\State \quad \textbf{if} $t < \tswitch$:
\State \qquad $\powercoeff{t}^2 = \powercoeff{lb}^2$
\State \quad \textbf{else}
\State \qquad $\powercoeff{t}^2 = \frac{\contr^{T-1-t}}{\sum_{\ell = \tswitch}^{T-1} \contr^{T-\ell - 1}}(T\Pavg - \tswitch\powercoeff{lb}^2)$ 
\State \textbf{end for} 
\State \Return $\{\powercoeff{t}\}_{t=0}^{T-1}$
\end{algorithmic}
\label{alg:poweralloc}
\end{algorithm}
}

In the following theorem, we show that using an exponentially increasing power allocation scheme, defined through Allocation Scheme~\ref{alg:poweralloc}, results in increased accuracy of the last iterate expected error. The power allocation scheme is obtained by minimizing the upper bound on the expected last iterate error with respect to $\{\powercoeff{t}\}$. 
\begin{theorem} \label{thm:PL} 
Consider Problem~$\cP_1$. 
Consider ${ \eta \in \paran{0,\frac{1}{2L}} }$ and  $\powercoeff{t}$ as:
$
    \powercoeff{t}^2 
    = \frac{\contr^{T-1-t}}{\sum_{k=0}^{T-1}\contr^{k}} T \Pavg
    \qquad \forall t \in [T-1]_0,
$
where $\contr := \sqrt{1-\mu\eta}$. 
If $\{n_t\}$ satisfy Assumption~\ref{asmp:noise} then \texttt{PAGD} with \emph{optimal allocation} ensure the following bound:
\begin{align}
    &\E[f(x_T)- \fopt] \cr 
    &\leq (1-\mu \eta)^T (f(x_0)- \fopt) + 
    \frac{\paran{\sum_{k=0}^{T-1} \contr^k }^2}{T} \frac{L G^2 \eta^2 \sigma_N^2}{\Pavg} \cr
    &\leq (1-\mu \eta)^T (f(x_0)- \fopt) +  \frac{4}{T}\frac{L G^2 }{\mu^2}\frac{ \sigma_N^2}{\Pavg}.
\end{align}
\end{theorem}
The proof for Theorem~\ref{thm:PL} is provided in Appendix~\ref{apndx:thm2}. The power allocation scheme leverages contraction of the error in each step using less power in the initial time steps. The resultant higher noise in the initial iterations is taken care of by the contraction at each step of the algorithm. 

Given the exponentially increasing power allocation, it is important to notice that the initial power allocation, ${\powercoeff{0}^2 = \frac{\contr^{T-1}}{\sum_{k=0}^{T-1}\contr^k}T\Pavg \leq \contr^{T-1}T\Pavg}$, decreases exponentially with the time horizon $T$. 
The effective noise added to the gradient in the update, $\frac{Gn_t}{\powercoeff{t}}$, is inversely proportional to the power coefficient $\powercoeff{t}$. 
Consequently, lower power allocation leads to higher variance in the noise being added to the gradients, which results in higher fluctuations in the update variables in the initial iterations of the algorithm. The increase in variance of the effective noise is relevant when we want to ensure the estimates stay within a compact set. 
\subsection{Optimization Problem: Locally Constrained}\label{subsec:opt_pow_all2}

Next we discuss the optimal power allocation and corresponding bound to functions satisfying local properties which will apply to the LQR problem. 
The effect of previously discussed increased variance is important when instead of satisfying $L$-smoothness throughout their domain, functions satisfy local $(L,D)$-smoothness as defined below. 
\begin{definition}[{Local Smoothness}] \label{def:local smooth}
A function $f:\R^d\to\R$ is said to be locally $(L,D)$-smooth over $\mathcal{X}\subseteq\mathbb{R}^d$ if $\|\nabla f(x)-\nabla f(y)\|_2 \le L\|x-y\|_2$ for all $x\in\mathcal{X}$ and all $y\in\mathbb{R}^d$ with $\|y-x\|_2 \le D$.
\end{definition}
Additionally, in the context of functions satisfying desired properties locally (local $(L,D)$-smoothness and $\mu$-PL condition within a compact set), we assume almost-sure boundedness constraint on noise, Assumption~\ref{asmp:noise_bounded}. This assumption is made to guarantee that the gradient descent paths remain within the intended range of values.

A natural way to deal with increased effective variance due to low power allocation is to set a lower bound, say $\powercoeff{lb}$, at every instant $t$, i.e., $\powercoeff{t}\geq \powercoeff{lb}$ for all $t \in [T]$. 
In the following theorem, we identify the optimal power allocation by minimizing the upper bound on the expected last iterate error with the additional lower bound constraints on power allocation. We provide a sufficient lower bound on the power allocation that ensures that the estimates at every time instant are within the desired set to enable the use of local $(L,D)$-smoothness property. Since we want to ensure that the estimates in every sample path of the random process are within the desired set, we have to use the almost-sure bounded property of the noise process. 
\av{\begin{theorem}[\texttt{PAGD} under Local Conditions] \label{thm:local_conv}
Consider ${f:\mathbb{R}^d\to\mathbb{R}_{\ge0}}$ and 
${\mathcal{X}=\{x\in\mathbb{R}^d:f(x)\le v\}}$, where $v\in\mathbb{R}_{\ge0}$. Suppose $f(\cdot)$ is $(L,D)$-smooth, satisfies $\mu$-PL condition over $\cX$, and has bounded gradients, $\norm{\nabla f(x)}_2 \le G$, for all $x\in\mathcal{X}$ and the random process $\{n_t\}$ satisfies Assumptions~\ref{asmp:noise} and~\ref{asmp:noise_bounded}.
\\
For a positive constant $\eta < \min\{D/G, 1/4L\} $, define $\contr = \sqrt{1-\mu\eta}$ and ${ \powercoeff{lb}^2:= G^2\Nmax^2 \max \paran{ \frac{\eta^2}{(D-G\eta)^2}, \frac{2}{\mu v}} }$. Suppose the average power budget satisfies $\Pavg \geq \powercoeff{lb}^2$. \\
Consider \texttt{PAGD} initialized with $x_0 \in \R^d$ such that $f(x_0)\le v/2$ and run with power allocation
\begin{align*}
    \{\powercoeff{t}\} =\texttt{CtG}(T,\Pavg,\mu,\eta,\powercoeff{lb})
\end{align*}
 as described in Allocation Scheme~\ref{alg:poweralloc}.
Then, for all $t\ge0$, $x_t\in\mathcal{X}$ and the expected error is bounded as follow: 
\begin{align}
    &\E[f(x_T)-\fopt] 
    \leq (1 - \mu\eta)^T (f(x_0)-\fopt) +\cr
   &\qquad  LG^2\eta^2 \sigma_N^2
    \paran{
    \frac{\sum_{t=0}^{\tswitch-1} \contr^{2(T-t-1)}}{\powercoeff{lb}^2} + \frac{{\sum_{t=\tswitch}^{T-1} \contr^{T-t-1}}}{T\Pavg - \tswitch\powercoeff{lb}^2}
    }. \nonumber
\end{align}
\end{theorem}
}

The proof of Theorem~\ref{thm:local_conv} is provided in Appendix~\ref{apndx:proof_local_conv}.

\textbf{Discussion:}
The power allocation scheme in Theorem~\ref{thm:local_conv} involves having constant power for $\tswitch$ iterations followed by the exponentially increasing scheme from Theorem~\ref{thm:PL}. For $t \geq \tswitch$, the power allocation can be expressed as  
$\powercoeff{k+\tswitch} = \frac{\contr^{T-\tswitch-1-k}}{\sum_{\ell=0}^{T-\tswitch-1} \contr^\ell} (T\Pavg - \tswitch\powercoeff{lb}^2)$
for ${k \in [T-\tswitch-1]_0}$. It is evident that for $t\geq \tswitch$ the power allocation is equivalent to the scheme of Theorem~\ref{thm:PL} with a time horizon $T-\tswitch$ and an average power budget $\frac{T\Pavg - \tswitch\powercoeff{lb}^2}{T-\tswitch}$ which is greater than $ \Pavg$ as long as $\tswitch < T$.   

$\powercoeff{lb}$ is the maximum of two terms- $(i)$ the lower bound $\powercoeff{lb}^2 \geq G^2\Nmax^2 \eta^2/(D-G\eta)^2$, which ensures that the difference in the estimates stays bounded that is  $\|x_{t+1}-x_t \| \leq D$  and $(ii)$ the lower bound, $\powercoeff{lb}^2 \geq 2G^2\Nmax^2/(\mu v) $, which ensures that the estimates stay in the desired sublevel set, $x_t\in \cX$ for all $t\in [T]$. 

For $(i)$, by selecting $\eta < \frac{D}{G\paran{1+\frac{\Nmax}{\Pavgcoeff}}} $ we can ensure that if $\powercoeff{lb}= \frac{G\Nmax\eta}{D-G\eta}$ then $\Pavg \geq \powercoeff{lb}$. 
For $(ii)$, depending on the specifics of the problem, we can adjust the lower bound while ensuring the result holds as long as the trajectory of the updates stays in the sublevel set $\cX$.

Note that one can obtain the result for using constant power allocation or increasing power scheme throughout by setting $\tswitch=T-1$ and $\tswitch=0$ respectively.  

\section{Solution to LQR problem}
For the LQR problem, the feasible set comprises the set of stabilizing controllers which lies in a sublevel set of the cost function $J(K)$.  

Next, we outline the key properties of the Linear Quadratic Regulator (LQR) problem that are relevant for our analysis. These properties are summarized from the results in \cite{fazel, cassel2021online}. To simplify notation, note the following definitions:
\begin{align}
&\beta_0 I\preceq R\preceq \beta_1 I,\ \beta_0 I\preceq Q\preceq \beta_1 I,\ \Sigma_w\succeq \sigma^2_w I,\cr 
&\norm{B}\le\psi,\ J(K^{\star})\le \frac{J}{4},\label{eqn:parameters in J(K)}
\end{align}
where $\beta_0,\beta_1,\sigma_w,J\in\mathbb{R}_{>0}$, $\psi\in\mathbb{R}_{\ge1}$ and $K^{\star}$ is the optimal solution to problem~\eqref{eqn:LQR obj J(K)}. Moreover, we assume without loss of generality that $\beta_1\le1$ (since one may always scale the cost matrices $Q,R$ by a positive real number). In addition, we construct a set 
\begin{equation}
\mathcal{K}=\{K\in\mathbb{R}^{m\times n}:J(K)\le J\},
\end{equation}
and impose $\mathcal{K}$ as the feasible set of $J(\cdot)$. 
In the following lemma, we characterize the properties of $J(\cdot)$ in terms of the parameters in Eq.~\eqref{eqn:parameters in J(K)}.  
\begin{lemma}{\cite[Lemma 5.1]{ye2024modelfreelearninglinearquadratic}}\label{lemma:properties of J(K)} The objective $J(\cdot)$ in problem~\eqref{eqn:LQR obj J(K)} satisfies: \\
\begin{enumerate}[a.]
\item \cite[Lemma~41]{cassel2020logarithmic} For any $K\in\mathcal{K}$, it holds that $\Vert A+BK \Vert^k\le\zeta(1-\xi)^k$ for all $k\in\mathbb{Z}_{\ge0}$ and $\Vert K \Vert \le\zeta$, where $\zeta\triangleq\sqrt{J/(\beta_0\sigma_w^2)}$ satisfies $\zeta\ge1$ and $\xi \triangleq1/(2\zeta^2)$.
\item \cite[Lemma~25]{fazel} For any $K\in\mathcal{K}$, it holds that $\norm{\nabla J(K)}_F\le G= \frac{2J}{\beta_0\sigma_w^2}\sqrt{(\sigma_w^2+\psi^2J)J}$.
\item\cite[Lemma~5]{cassel2021online} $J(\cdot)$ is $(L,D)$-locally smooth with $D=\frac{1}{\psi\zeta^3}$ and $L=112\sqrt{n}J\psi^2\zeta^8/\beta_0$, i.e., $\norm{\nabla J(K^{\prime})-\nabla J(K)}_F\le L\norm{K^{\prime}-K}_F$ for all $K\in\mathcal{K}$ and all $K^{\prime}\in\mathbb{R}^{m\times n}$ with $\norm{K^{\prime}-K}\le D$.
\item\cite[Lemma~11]{fazel} $J(\cdot)$ satisfies the gradient-domination property with $\mu=2J/\zeta^4$, i.e., $\Vert \nabla J(K)\Vert^2\ge 2\mu(J(K)-J(K^{\star}))$ for all $K\in\mathcal{K}$, where $K^{\star}=\argmin_{K\in\mathcal{K}}J(K)$.
\end{enumerate}
\end{lemma}

Using the above lemma to determine the properties for the problem parameters $\mu, L, D$, and $ G$ and applying Theorem~\ref{thm:local_conv}, we get the following theorem. 
\begin{theorem}\label{thm:lqr_thm} Consider $\cP_0$ for which the values of the parameters $\mu, L, G,$ and $D$ are as stated in Lemma~\ref{lemma:properties of J(K)}. Let the policy gradient method be initialized with $K_0$ and constant step-size $\eta$ satisfying $0<\eta < \min\{D/G, 1/4L\} $. 
Define $\contr = \sqrt{1-\mu\eta}$ and ${ \powercoeff{lb}^2:= G^2\Nmax^2 \max \paran{ \frac{\eta^2}{(D-G\eta)^2}, \frac{2}{\mu J}} }$. Suppose the average power budget satisfies $\Pavg \geq \powercoeff{lb}^2$. 
Using \texttt{PAGD} with \emph{optimal allocation}, Allocation Scheme~\ref{alg:poweralloc}, for the LQR problem results in 
    \begin{align*}
        &\E[J(K_t) - J^{\star}]
        \leq 
        (1 - \mu \eta)^T (J(K_0) - J^{\star})  +\cr
        &
         LG^2\eta^2 \sigma_N^2
    \paran{
    \frac{\sum_{t=0}^{\tswitch-1} \contr^{2(T-t-1)}}{\powercoeff{lb}^2} + \frac{{\sum_{t=\tswitch}^{T-1} \contr^{T-t-1}}}{T\Pavg - \tswitch\powercoeff{lb}^2}
    }.
    \end{align*}
\end{theorem}

\av{\section{Conclusion}
We studied the policy gradient method for the LQR problem when the gradients are transmitted by an agent with limited power budget over a noisy communication channel. To address this, we proposed closed-form power allocation strategies that follow a hybrid of constant and geometrically increasing structure. These allocations are derived within an optimization framework for two class of functions that satisfy smoothness and PL conditions either globally or locally. Rather than directly minimizing the expected suboptimality, we optimize a tractable upper bound on the expected error in the function value. Finally, we apply our approach to the policy gradient setting in LQR problem, demonstrating how the optimized allocation effectively supports convergence under communication constraints.
}
\bibliographystyle{ieeetr}
\bibliography{biblio}

\appendices
\section{Solution to the Optimization Problem}\label{apndx:OptProblem}
\begin{lemma}\label{lem:OptProblem2}
    Let  $\{a_i\} $ be a sequence of increasing positive constants for $ i \in [n]$. Consider the following optimization problem:
\begin{align}
    \min_{w_1, w_2, \dots, w_n} & \quad \sum_{i=1}^{n} \frac{a_i}{w_i} \label{eq:optproblem} \\
    \text{subject to:} & \quad \sum_{i=1}^{n} w_i \leq K, 
                        \quad w_i \geq C_L, \quad \forall i \in \{1, 2, \dots, n\}.\nonumber
\end{align}
The minimum value of the objective function is achieved at 
\begin{align}
    \sum_{i=1}^{i_{\cS}-1} \frac{a_i}{C_L} + \sum_{i=i_{\cS}}^n \sqrt{a_i\lambda} 
    = 
    \frac{\sum_{i=1}^{i_{\cS}-1} a_i}{C_L} 
    + \frac{\paran{\sum_{i=i_{\cS}}^n \sqrt{a_i}}^2}{K-(i_{\cS}-1)C_L} \nonumber
\end{align}
with the optimal solution given by 
\begin{align}
    w_i = 
    \begin{cases}
        \sqrt{\frac{a_i}{\lambda}} & \text{ if }i \geq i_{\cS} \\
         C_L & \text{ if } i < i_{\cS}
    \end{cases}, \qquad \forall i\in [n],
    \end{align}
where $ i_{\cS}$ is defined as $i_{\cS}  = \min\{i\in[n] \mid \sqrt{\frac{a_i}{\lambda(i)}} \geq C_L\}$ for  $\lambda(j) = \left(
    \frac{\sum_{i=j}^n \sqrt{a_i} }{K- (j-1)C_L }
    \right)^2.
$
\end{lemma}
{Note that $i_{\cS}$ does not have a closed form expression and needs to be determined by performing a (binary) search. $i_{\cS}$ is the index up to which we assign the variables $w_i$ with the minimum threshold to satisfy the constraint $w_i \geq C_L$ and after which (possibly) varying value of $w_i$ comes into play. We need ${nC_L \leq K}$ for the feasibility of $w_i\geq C_L$ to hold for all $i\in [n]$. If $K \geq nC_L$, then it is easy to see that $i_{\cS}\leq n$.} 

\begin{proof}
The optimization problem is convex. To solve the optimization problem, we begin with the Lagrangian with $\lambda>0, \mu_i>0$ for $i\in [n]$:
\begin{align}
\cL(w, \lambda, \mu) = \sum_{i=1}^n \frac{a_i}{w_i} + \lambda \left( \sum_{i=1}^n w_i - K \right) - \sum_{i=1}^n \mu_i (w_i - C_L). \nonumber
\end{align}
The KKT conditions for the problem give us the following. 
\begin{enumerate}
\item Derivative with respect to $w_i$:
$\del{\mathcal{L}}{w_i} = -\frac{a_i}{w_i^2} + \lambda - \mu_i = 0.$
Solving for $w_i$, we get:
   \begin{align*}
    w_i = \sqrt{\frac{a_i}{\lambda - \mu_i}}, \quad \text{if } \lambda - \mu_i > 0.    
   \end{align*}
\item Complementary slackness for the constraints: ${\mu_i (w_i - C_L) = 0 \quad \forall i\in [n].}$
    If $\mu_i > 0$, then $w_i = C_L$. Otherwise, $w_i = \sqrt{\frac{a_i}{\lambda}}$.
    \item The budget constraint:
   $
   \sum_{i=1}^n w_i \leq K.
   $
\end{enumerate}

The solution depends on whether the unconstrained optimal satisfies the desired constraint, i.e., $\sqrt{\frac{a_i}{\lambda}} \geq C_L$. If 
$\sqrt{\frac{a_i}{\lambda}} < C_L$, then $w_i = C_L$.

Define the set of indices where the unconstrained solution, $\sqrt{\frac{a_i}{\lambda}}$, is greater than or equal to $C_L$:
\begin{align}
      \cS :=
      \left\{i\in[n] \bigg | \sqrt{\frac{a_i}{\lambda}} \geq C_L\right\}. \nonumber
\end{align}
 Since $\{a_i\}$ is a sequence of increasing positive constants, define 
$  i_{\cS} := \min\{i\in[n] \mid \sqrt{\frac{a_i}{\lambda}} \geq C_L\}, $
with the convention being $i_{\cS} := n+1$ if $\sqrt{\frac{a_n}{\lambda}} < C_L$. Then ${\cS = \{i\in[n] | i\geq i_{\cS}\}.}$

For $i\in [n]$ set $w_i$ as follow:
\begin{align}
    w_i = 
    \begin{cases}
        \sqrt{\frac{a_i}{\lambda}} & \text{ if }i \geq i_{\cS}, \\
         C_L & \text{ if } i < i_{\cS}
    \end{cases}.\nonumber
\end{align}

Using the above assignment the budget constraint implies 
   $
   \sum_{i =i_{\cS}}^n \sqrt{\frac{a_i}{\lambda}} + \sum_{i =1}^{i_{\cS}-1} C_L = K.
   $ 
   Solving for $\lambda$ we get 
$    \lambda = \left(
    \frac{\sum_{i=i_{\cS}}^n \sqrt{a_i} }{K- (i_{\cS}-1)C_L }
    \right)^2.
$
\end{proof}
The following corollary follows from Lemma~\ref{lem:OptProblem2}.
\begin{corollary}\label{lem:OptProblem}
    Let  $a_i > 0 $ be constants for $ i \in [n]$, and consider the optimization problem~\ref{eq:optproblem} with $C_L = 0$. 
    The minimum value of the objective function is achieved at $ \frac{\paran{ \sum_{i=1}^n \sqrt{a_i}}^2}{K}$, with the optimal solution given by ${w_i = K\frac{\sqrt{a_i}}{\sum_{k=1}^n \sqrt{a_k}} }$ for every $i \in [n] $.
\end{corollary}
\section{Proof of Theorem~\ref{thm:PL}}\label{apndx:thm2}
\begin{proof}
Recall that $L$-smoothness implies that for all $x,y \in \R^d$, 
$ f(y) \leq f(x) + \inp{\nabla f(x)}{y-x} + \frac{L}{2}\norm{y-x}^2. $
For $t\in [T-1]_0$, the function's value at the estimate can be upper-bounded as
\begin{align}
&f(x_{t+1})  \overset{(a)}= f(x_t - \eta \dec(g_t))
\cr &\overset{(b)}
\leq f(x_t) - \eta\inp{\nabla f(x_t)}{\dec(g_t)} + \frac{L}{2} \eta^2\| \dec(g_t)\|^2 \cr
&= f(x_t) - \eta \|\nabla f(x_t)\|^2 - \frac{G\eta}{\powercoeff{t}}\inp{\nabla f(x_t)}{n_t} 
\cr &\qquad + \frac{L \eta^2}{2}\norm{\nabla f(x_t) + \frac{G}{\powercoeff{t}}n_t}^2 \cr
&\overset{(c)}\leq 
f(x_t) - \eta \|\nabla f(x_t)\|^2 
- \frac{G\eta}{\powercoeff{t}}\inp{\nabla f(x_t)}{n_t} 
\cr &\qquad + L\eta^2\norm{\nabla f(x_t)}^2 + \frac{LG^2\eta^2}{\powercoeff{t}^2}\norm{n_t^2}, \label{ineq:innerp}
\end{align}
where $(a)$ follows from eq.~\eqref{eq:updateRule} $(b)$ follows from the $L$-smoothness, and $(c)$ follows from the inequality $\|a+b\|^2 \leq 2\|a\|^2 +2 \|b\|^2$ which holds for any $a,b \in \R^d$.  

Taking expectation conditioned on sigma-field based on the information till time $t$, 
$\cF_t=\sigma(\bigcup_{k=0}^{t-1} \{x_k,n_k\}\cup \{x_t\})$.
and using the stochastic properties of the noise specified in Assumption~\ref{asmp:noise} implies $\E[n_t \mid \cF_t] = \vec{0}$ and $\E[\|n_t\|^2 \mid \cF_t] = \sigma_N^2$. Taking conditional expectation on \eqref{ineq:innerp} we get 
\begin{align}
    \E[f(x_{t+1}) | \cF_t] 
    &
    \leq f(x_t) - \paran{\eta -L\eta^2 } \|\nabla f(x_t)\|^2 
    + \frac{LG^2\eta^2}{\powercoeff{t}^2}\sigma_N^2. \label{ineq:Ef_recursion}
\end{align} 

Define the distance of the function's value at the current estimate from the optimal as ${z_t := f(x_t) - \fopt}$. 
Using the gradient-dominance property, $\norm{\nabla f(x)}^2 \geq 2\mu (f(x) - \fopt) $, we derive 
\begin{align}
    \E[z_{t+1} | \cF_t] 
    &\leq \paran{1- 2\mu\eta + 2\mu L\eta^2} z_t + \frac{LG^2\eta^2}{\powercoeff{t}^2}\sigma_N^2 \cr 
    &=: (1-b)z_t + \frac{c}{\snr{t}},\label{ineq:Ez_recursion}
\end{align}
where $b := 2\mu \eta - 2{\mu L \eta^2}$, $c := {L G^2 \eta^2}$, and $\snr{t} := \frac{\powercoeff{t}^2}{\sigma_N^2}$ is the Signal-to-Noise power Ratio (SNR) for transmission at time $t$. 
Unrolling \eqref{ineq:Ez_recursion} at time $T$, we get 
\begin{align}
    \E[z_T] \leq (1-b)^T z_0  + \sum_{t=0}^{T-1} (1-b)^{T-1-t} \frac{c}{\snr{t}}. \label{ineq:unrolled1}
\end{align}

\textbf{Optimization problem for $\{\snr{t}\}$:}
\begin{align}
    \min_{\snr{0},\snr{1},\dots,\snr{T-1}} & \quad \sum_{t=0}^{T-1} (1-b)^{T-t-1} \frac{c}{\snr{t}} \cr 
    \text{subject to} & \quad \frac{\sum_{t=0}^{T-1} {\snr{t}}}{T}\leq 
    \bar{\snr{}};
    \quad {\snr{t}}\geq 
    0, \forall t\in [T-1]_0, \nonumber
\end{align}
where $\bar{\snr{}} := \frac{\Pavg}{\sigma_N^2}$. 
According to Corollary~\ref{lem:OptProblem} the optimal selection of $\{\snr{t}\}$, for ${t \in [T-1]_0}$, is given by 
\begin{align}
    \snr{t} &= T\bar{\snr{}}\frac{(1-b)^{\frac{T-t-1}{2}}}{\sum_{t=0}^{T-1}(1-b)^{\frac{T-t-1}{2}}}\label{eq:snr_opt} 
    = \frac{(1-\sqrt{1-b})(1-b)^{\frac{T-t-1}{2}}}{1- (\sqrt{1-b})^{\frac{T-1}{2}}}T\bar{\snr{}} 
\end{align}
The minimized error term is \red{$c
\frac{(1-(\sqrt{1-b})^{T})^2}{T(1-\sqrt{1-b})^2} $}. 

To ensure contraction of the initial error we need the step-size $\eta$ to satisfy $\eta L \leq 1$ which guarantees that $|1-b|<1$. Note when $\eta = \frac{1}{2L}$ the parameter $b= \frac{\mu}{2L}$ is maximized.
Finally note that if $\eta L \leq \frac{1}{2}$, then $b = 2\mu \eta (1-{L\eta})\geq \mu \eta$ and thus giving $(1-b) < (1-\mu\eta)$.

Using $\eta<1/(2L)$ and substituting back $c = {L G^2 \eta^2}$, we obtain the following bound on the expected error when the power allocation is according to \eqref{eq:snr_opt}.
\begin{align}
    &\E[z_T] 
    \leq (1-b)^T z_0 + \frac{(1-(\sqrt{1-b})^{T})^2}{(1-\sqrt{1-b})^2}\frac{L G^2 \eta^2 \sigma_N^2}{T \Pavg} \cr 
    &\leq (1-\mu\eta)^T z_0 + \frac{(1-(\sqrt{1-\mu\eta})^{T})^2}{(1-\sqrt{1-\mu\eta})^2}\frac{L G^2 \eta^2 \sigma_N^2}{T \Pavg}. \label{ineq:zt1}
\end{align}
Finally using the fact that $\frac{1 - a^T}{1-a}\leq \frac{1}{1-a}$ and ${\sqrt{1-a}\leq 1-\frac{a}{2}}$ we get the bound $\frac{(1-(\sqrt{1-\mu\eta})^{T})^2}{(1-\sqrt{1-\mu\eta})^2} \leq \frac{4}{\mu^2\eta^2}$. Utilizing the bound in inequality~\ref{ineq:zt1} we get 
\begin{align*}
    \E[z_T] \leq (1-\mu\eta)^T z_0 + \frac{4}{T}\frac{LG^2}{\mu^2} \frac{\sigma_N^2}{\Pavg}.
\end{align*}
\end{proof}
\section{Proof of Theorem~\ref{thm:local_conv}}\label{apndx:proof_local_conv}
\begin{proof} 
First, we establish that $\|x_{t+1}- x_t\|\leq D$ for all $t\in [T-1]_0$. We know that 
\begin{align}
\|x_{t+1} - x_t\| 
&= \|\eta \dec(g_t)\| \leq \eta\|\nabla f(x_t)\| + \frac{G\eta}{\powercoeff{t}} \|n_t\|  \cr 
&\leq \eta G + \eta\frac{G\Nmax}{\powercoeff{t}} 
\leq G\eta\paran{ 1 + \frac{\Nmax}{\powercoeff{t}}} 
\leq G\eta\paran{ 1 + \frac{\Nmax}{\powercoeff{lb}}} \cr
&\leq D, \label{eq:Dbound}
\end{align}
where we use the fact that $\powercoeff{lb} \geq G\Nmax \frac{\eta}{D-G\eta} = \frac{\Nmax}{\frac{D}{G\eta}-1}$. 
Next, we establish that the sequence of estimates $\{x_t\}$ lies in the desired sublevel set $\cX$, i.e., $f(x_t) \leq v$ for all ${t\in [T-1]_0}$. 
From eq.~\eqref{ineq:innerp} we know 
\begin{align}
    z_{t+1}
    &\leq
    z_t - \eta \norm{\nabla f(x_t)}^2 
    - \frac{G\eta}{\powercoeff{t}}\inp{\nabla f(x_t)}{{n_t}}
    \cr & \quad + {L\eta^2} \norm{\nabla f(x_t)}^2 + \frac{LG^2\eta^2}{\powercoeff{t}^2}\norm{n_t}^2 \cr 
    &\leq
    z_t - \eta \norm{\nabla f(x_t)}^2+ {L\eta^2} \norm{\nabla f(x_t)}^2 
     \cr &\quad + \frac{G\eta}{2}\paran{\frac{\norm{\nabla f(x_t)}^2}{G} + G\frac{\norm{n_t}^2}{{\powercoeff{t}}^2}}+ \frac{LG^2\eta^2}{\powercoeff{t}^2}\norm{n_t}^2 \cr
     &\leq 
     \paran{1-2\mu\eta \paran{\frac{1}{2}-L\eta}}z_t +G^2\eta\paran{\frac{1}{2}+{L\eta}}\frac{\norm{n_t}^2}{\powercoeff{t}^2}.\nonumber
\end{align}
If \red{$\eta L \leq \frac{1}{4}$}, we can upper bound the expected error as follow
\begin{align}
    z_{t+1}
    &\leq 
     \paran{1-\frac{\mu\eta}{2}}z_t +G^2\eta\frac{\norm{n_t}^2}{\powercoeff{t}^2} \leq 
     \paran{1-\frac{\mu\eta}{2}}z_t +G^2\eta\frac{\Nmax^2}{\Cmin^2}. \nonumber
\end{align}
\red{
We now prove $z_{t}\leq v$ for all $t \in [T]_0$ by induction. The base case $z_0 \leq v$ holds true. 
Assume the induction hypothesis, $z_t \leq v$ for some $t\in [T-1]_0$. Since $\powercoeff{lb}^2 \geq \frac{2G^2\Nmax^2 }{v\mu}$ we have $ z_{t+1}\leq v$. Therefore by induction we know $z_t \leq v$, and thus $f(x_t) \leq v$, for all $t\in [T]_0$.
}

Since the $(L,D)$-smoothness and $\mu$-PL condition are satisfied for all ${t \in [T-1]_0}$ the inequality \eqref{ineq:unrolled1} holds:
\begin{align}
    \E[z_T] \leq (1-b)^T z_0  + \sum_{t=0}^{T-1} (1-b)^{T-t-1} \frac{c}{\snr{t}}. \nonumber
\end{align}
The sequence of $\powercoeff{t}$ can be chosen to solve the following optimization problem:
\begin{align}
    \min & \quad \sum_{t=0}^{T-1} (1-b)^{T-t-1} \frac{c}{\snr{t}} \cr 
    \text{s.t.} &\quad \frac{\sum_{t=0}^{T-1} {\snr{t}}}{T}\leq 
    \bar{\snr{}}; \quad \snr{t} \geq \snr{lb}  \qquad \forall t \in [T-1]_0, \nonumber
\end{align}
where 
\red{$\snr{lb}:= 
\frac{\Cmin^2}{\sigma_N^2}
=G^2 \frac{\Nmax^2}{\sigma_N^2}
\max \paran{
\frac{\eta^2}{(D-G\eta)^2},
\frac{2}{\mu v}
}$}, 
${1-b=1- 2\mu\eta(1 - {L\eta})\leq 1-\frac{3}{2}\mu\eta\leq 1-\mu\eta}$, and ${c = LG^2\eta^2}$. 
The solution of the optimization problem from Lemma~\ref{lem:OptProblem2} gives
$\snr{t} = \snr{lb}$ for $t < \tswitch$ and 
\begin{align}
    \! \snr{t} \!=\!
        \sqrt{\frac{(1-b)^{T-t-1}}{\lambda(\tswitch)}}
        \!=\! \frac{\sqrt{(1-b)^{T-t-1}}}{\sum_{\ell=\tswitch}^{T-1} \sqrt{1-b}^{T-\ell-1}} (T\bar{\snr{}} - \tswitch\snr{lb}), \nonumber 
\end{align}
 if $t \geq \tswitch$, 
where $\lambda(\tswitch) := \paran{ \frac{\sum_{t=\tswitch}^{T-1} \sqrt{(1-b)^{T-t-1}}}{T\bar{\snr{}}-(\tswitch)\snr{lb}} }^2$ with 
${\tswitch = \min\{t\in [T-1]_0 \mid \sqrt{\frac{(1-b)^{T-t-1}}{\lambda(t)}}\geq \snr{lb}\}.}$
Using this power allocation, and $\contr =\sqrt{1-\mu\eta}$, the expected distance from the optimal is 
\begin{align}
    &\E[f(x_T)-\fopt] 
    \leq (1 - \mu\eta)^T (f(x_0)-\fopt) \cr
    &\qquad + LG^2\eta^2 
    \left(
    \frac{\sum_{t=0}^{\tswitch-1} \contr^{2(T-t-1)}}{\snr{lb}} +
    \frac{\sum_{t=\tswitch}^{T-1} \contr^{T-t-1}}{T\bar{\snr{}} - \tswitch\snr{lb}}
    \right). \nonumber
\end{align}
\end{proof}
\end{document}